\input amstex.tex
\documentstyle{amsppt}
\input graphicx.tex 

\magnification=1200
\hsize=150truemm
\vsize=224.4truemm
\hoffset=4.8truemm
\voffset=12truemm

\TagsOnRight
\NoBlackBoxes
\NoRunningHeads

\def\Square{\rlap{$\sqcup$}$\sqcap$}
\def\cqfd {\quad \hglue
7pt\par\vskip-\baselineskip\vskip-\parskip
{\rightline{\Square}}}

\define\R{{\bold R}}
\define\Z{{\bold Z}}
\define\mi{^{-1}}
\let\thm\proclaim
\let\fthm\endproclaim
\let\inc\subset 
\let\ev\emptyset
\let\ov\overline

\newcount\refno
\global\refno=0
\def\nextref#1{\global\advance\refno by
1\xdef#1{\the\refno}}
\def\bref {\ref\global\advance\refno by 1\key{\the\refno}}



\nextref\Bow
\nextref\FoGBS
\nextref\Fo
\nextref\FoJSJ
\nextref\GHMR
\nextref\Gui
\nextref\GuiLe
\nextref\LeOuthyp
\nextref\LeGBS
\nextref\Wh

\topmatter

\abstract We extend Forester's rigidity theorem so as to give
a complete characterization of rigid group actions on trees (an
action is rigid if it is the only reduced action   in its
deformation space, in particular it is invariant under
automorphisms preserving the set of elliptic subgroups). 
\endabstract

\title Characterizing  rigid simplicial actions on trees
      \endtitle

\author  Gilbert Levitt
\endauthor

\endtopmatter

\document

Let $T$ be a simplicial tree with a cocompact action of a group
$G$ (i.e. the Bass-Serre tree associated to a  decomposition
    of $G$ as a finite graph of groups). A $G$-tree $T'$ is a {\it
deformation\/} of $T$ if it may be obtained from $T$ by a finite
sequence of expansions and collapses (elementary moves
coming from the canonical isomorphism
$A*_BB\simeq A$). These moves do not change the set of
elliptic subgroups (a subgroup is {\it elliptic\/} if it fixes a
point in the tree). Conversely, Forester proved [\Fo] that any
cocompact $G$-tree $T'$ with the same elliptic  subgroups as
$T$ is a deformation of
$T$.

Since an expansion makes the tree more complicated, and a
collapse makes it simpler, it is natural to restrict to reduced
trees. A tree $T$ is {\it reduced\/} if one cannot perform a
collapse on $T$. Equivalently,
$T$ is reduced if, whenever an edge $e$ has the same stabilizer
as one of its endpoints, then both endpoints of
$e$ are in the same $G$-orbit (i.e\. $e$ projects onto a   loop in
the quotient graph).

The  tree $T$ is {\it rigid\/} if it is reduced, and it is the only
reduced tree in its deformation space (up to equivariant
isomorphism). In other words, all deformations of
$T$ (trees with the same elliptic subgroups as $T$) may be
reduced to $T$ by collapse moves. Rigidity provides a
canonical element
$T_{red}$ in the deformation space; in particular, any
automorphism of $G$ that preserves the set of elliptic
subgroups leaves $T_{red}$ invariant (see [\Bow, \FoJSJ,
\GHMR, 
\GuiLe, \LeOuthyp, \LeGBS] for examples and applications to  
JSJ splittings and automorphisms).

Forester proved   that ``strongly-slide-free'' trees are rigid
([\Fo], see also [\Gui]). The purpose of this note is to extend
Forester's theorem. Our extension is optimal: we obtain  a
complete characterization of rigid trees. 

Before stating our result, let us illustrate it  on generalized
Baumslag-Solitar groups [\FoGBS, \LeGBS]. Note that these groups
have been classified up to quasi-isometry [\Wh]. The rigidity we are
studying here is not quasi-isometric rigidity of groups, as the group
is fixed once and for all.

We consider a finite graph of groups with each vertex and edge
group isomorphic to $\Z$. It is pictured as a labelled graph,
each label being the index of the edge group in the vertex
group (see figure 1). [One should  allow negative labels, but we
will not bother here.]  

\midinsert
\centerline 
{\includegraphics[scale=.9]
{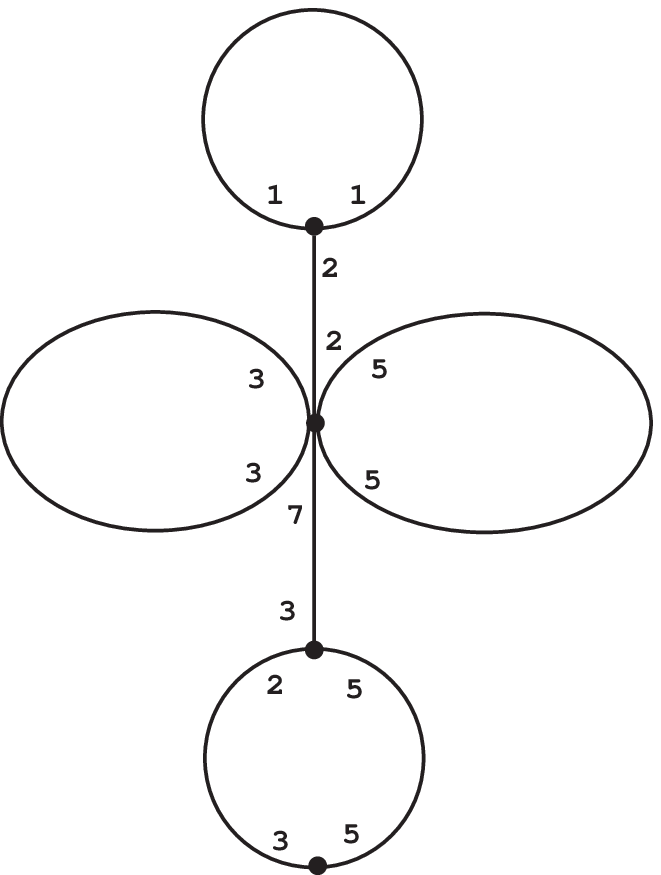}}
\botcaption
 {Figure 1}{A graph of groups whose 
universal cover  is a rigid tree}
\endcaption
\endinsert


The associated Bass-Serre tree $T$ is reduced if and only if the
label $1$ appears only on   loops (edges with both endpoints
equal). It is strongly slide-free (in the sense of [\Fo]) if and
only if there is no divisibility relation: a divisibility relation is
a relation $p\mid q$, where $p,q$ are labels at the same vertex.

We will show that the tree   associated to a labelled graph is
rigid if and only if the graph
   is as on figure 1. Namely, a  divisibility relation $p\mid q$ at a
vertex $x$ is allowed if $p=q$ and $p,q$ are carried by one
loop. If $x$ has valence $3$, with a loop labelled
$(1,1)$ and a  third label
$r$, we allow the relations $1\mid r$. No other divisibility
relation is allowed.

There is one exception to the statement just given. The tree
associated to the standard presentation of the solvable
Baumslag-Solitar group $BS(1,s)$ as an HNN-extension   is rigid
if and only if $s=1$ or $s$ is prime.

We shall now give the general statement.  Let $T$ be a
$G$-tree. If $e$ is an oriented edge, $\bar e$ denotes $e $ with
the opposite orientation. We always assume that $G$ acts
without inversions. We denote by
$G_v$ (resp\. $G_e$) the stabilizer of a vertex $v$ (resp\. an
edge $e$).

A   $G$-tree $T$ is (associated to) {\it an ascending
HNN-extension\/} if the quotient graph of groups has one
vertex and one edge, and  furthermore at least one of the two
monomorphisms from the edge group to the vertex group is
onto. If $T$ is reduced, this is equivalent to the existence of a
$G$-fixed end in $T$, and also to the length function being the
absolute value of a homomorphism $G\to\R$.

\thm{Theorem 1 (general case)} Let $T$ be a reduced 
cocompact
$G$-tree which is not an ascending HNN-extension. It   is rigid if
and only if, given two oriented edges
$e, f$ with the same   origin $v$ such that $G_e\inc G_f$, one of
the following holds:
\roster
\item
$e$ and $f$ are in the same $G_v$ orbit.
\item
$e$ and $\bar f$ are in the same $G$-orbit, and $G_e=G_f$.
\item there is an edge $f'$ with origin $v$,  in the same
$G$-orbit as $\bar f$, such that
$G_f=G_{f'}=G_v$. Furthermore, there are only three
$G_v$-orbits of edges   with origin
$v$ (those of
$e,f,f'$).
\endroster
\fthm

Strongly slide-free trees are those for which only $(1)$ occurs.

\thm{Theorem 2 (exceptional case)} Let $T$ be the Bass-Serre
tree of  an ascending HNN-extension
$G =\langle A,t\mid tat\mi=\varphi (a)\rangle
$, with $\varphi :A\to A$ an injective homomorphism.

The tree $T$ is rigid if and only if, for every subgroup $H\inc
A$  containing $\varphi (A)$, there exist $i,j\ge0$ and $a_0,
a_1\in A$ such that $a_0$ conjugates $\varphi ^i(A)$ to
$\varphi ^j(H)$, and 
$\varphi (a_0\mi)a_0\varphi ^{i+1}(a_1)$ centralizes  $\varphi
^{i+1}(A)$.
\fthm

If $\varphi (A)$ is either $A$ or a maximal proper subgroup of
$A$, the tree is rigid. The converse holds when $A$ is abelian.

Theorem 1 is proved in the first two sections.
Theorem 2 is proved in the third one.

\head Rigidity \endhead

We prove the ``if'' direction of Theorem 1.    Let
$T'$ be reduced, with the same elliptic subgroups as $T$. We
show
$T'=T$.

$\bullet$ First assume that case $(3)$ in the statement of Theorem 1
does not occur. Our arguments generalize those of [\Gui]. 

We define a map $f:T\to T'$ in the
following way. We choose a representative $v_i$ for each orbit of
vertices of $T$, and we let $f(v_i)$ be a vertex  of $T'$ fixed by
$G_v$.  We then extend  $f$ to a $G$-equivariant  map
  sending vertex to vertex, and linear on edges.

Since $T$ is as in Theorem 1, and case $(3)$ does not
occur,
  no vertex stabilizer $G_v$ of $T$ can fix
an edge. As in [\Gui, p\. 324], it follows that distinct vertices of
$T$ have distinct images in $T'$. Also note that the stabilizer of
$f(v)$ equals the stabilizer of $v$.

The key point is to show that there is no folding,
a fold being a pair of edges $e_1,e_2$ with origin
$v$ such that
$f(e_1)\cap f(e_2)$ is a non-degenerate segment.

Suppose there is a fold. 
  Let $w'$ be the vertex of $T'$ adjacent to $f(v)$ on $f(e_1)\cap
f(e_2)$. 
The argument  in the proof of Lemma 2.1 in [\Gui] shows that the
stabilizer $G_{w'}$ of $w'$ cannot be contained in $G_v$. Since  
$G_{w'}$ is elliptic, it fixes a vertex
$w\neq v$ in $T$. Let $e$ be the initial edge of
$vw$.  

Any element of $G_v$ fixing $w'$ also fixes $e$. 
In particular,  $G_{e }$ contains both $G_{e_1}$ and
$G_{e_2}$. Since $(1)$ or $(2)$ holds for each of the pairs
$(e_1,e )$ and $(e_2,e )$, we know that
$e_1,e_2,e $ are in the same orbit as non-oriented edges. In
particular, $f(e_1)$ and $f(e_2)$ have the same length.

Say that the fold between $f(e_1)$ and $f(e_2)$ is of type
$(1)$ if
      $e_1$, $e_2$ are in the same  $G_v$-orbit. In this case 
$G_{e_1}$ and $G_{e_2}$ are properly contained in $G_{e }$,
because elements of $G_v$ taking $e_1$ to
$e_2$ belong to $G_{e }$ (they fix $w'$). In particular, $e $ is in the
same
$G_v$-orbit as
$e_1$ and
$e_2$.

If the fold is of type $(2)$ (i.e\. $e_1$ and $\bar e_2$ are in the
same
$G$-orbit), we claim that $G_{e_1}=G_{e_2}$. We may assume
that $e_1$ and $e $ are in the same $G_v$-orbit. Since
$G_{e_2}\inc G_{e }$, we have $G_{e_2} = G_{e }$. But then
$G_{e_1}\inc G_{e_2}$, and finally $G_{e_1}= G_{e_2}$. Also
note that   $f(e_1)$ and  $f(e_2)$ are folded along strictly less
than half their (common) length: otherwise, an element of $G$
sending $e_1$ to
$\bar e_2$ would fix a point in $T'$ but not in $T$.

We can now extend Lemma 2.2 of [\Gui]:

\thm{Lemma}{\it If $e_1,e_2,e_3$ are three consecutive
(non-oriented)  edges in $T$, then $f(e_1)\cap f(e_2)\cap
f(e_3)=\ev$.}
\fthm

\demo{Proof} Denote
$v_1=e_1\cap e_2$ and $v_2=e_2\cap e_3$. There is a problem
only if there are folds both at $v_1$ and at $v_2$. If both folds
are of type $(1)$, the argument of [\Gui] applies. If both folds
are of type $(2)$, we simply use the fact that less than half is
folded. We complete the proof by ruling out the possibility of
mixed folds: type
$(1)$ at
$v_1$, type $(2)$ at $v_2$.

Orient $e_1,e_2,e_3$ so that $e_1$ and $e_2$ have origin
$v_1$, and $e_3$ has origin
$v_2$. They are in the same $G$-orbit as oriented edges. Let
$e_0$ be the image of
$e_2$ by a group element taking $e_3$ to $e_1$. Its terminal   
endpoint is $v_1$. Since the fold at
$v_2$ is of type $(2)$, we have $G_{e_2}=G_{e_3}$, and
therefore
$G_{e_0}=G_{e_1}$.

Now consider an  edge $e $ (with origin $v_1$) associated to
the fold
$e_1,e_2$ as above. It is in the $G_{v_1}$-orbit of $e_1$ (but
not of $\bar e_0$). Furthermore
$G_e$ properly contains $G_{e_1}$, hence also $G_{e_0}$. This 
shows that
$\bar e_0, e$ do not satisfy any of conditions $(1)$ or
$(2)$ of the theorem, a contradiction. \cqfd

\enddemo

Deducing $T'=T$ from this lemma  is  as in [\Gui, pp\.  326-327].

\smallskip
$\bullet$ We now have to allow  case $(3)$ in the statement of
Theorem 1 to occur. Let $e,f,f',v$ be as in  case $(3)$. Since $T$ is
reduced,
$G_e$ is properly contained in $G_v$. It follows that the
normalizer of $G_v$ is a semi-direct product
$N(G_v)=G_v\rtimes
\langle t \rangle$. It acts on a line $L_v\inc T$ (containing $f$
and $f'$), with $G_v$ acting as the identity and $t$ as a unit
translation (taking $\ov f$ to $f'$). Every point of $L_v$ has
stabilizer $G_v$.

We claim that the translation axis $L'_v$ of $t$ in $T'$ has
similar properties. First,
$G_v$ is the identity on $L'_v$ because $G_v$ is elliptic in $T'$
and
$t$ normalizes $G_v$. Since $G_v$ is a maximal elliptic
subgroup, points of $L'_v$ have stabilizer equal to
$G_v$. In particular, $gL'_v\cap L'_v=\ev$ if $g\notin N(G_v)$.
Adjacent vertices on
$L'_v$ are in the same
$G$-orbit because
$T'$ is reduced. They are in the same orbit under $t$ because
$N(G_v)$ is generated by $G_v$ and
$t$, so $t$ acts on $L'_v$ as a unit translation.

For each vertex $v$ as in case $(3)$ of Theorem 1 (there may
be several $G$-orbits of them), collapse $L_v$ to a point.
Similarly, collapse each $L'_v$ (noting that collapsed lines are
pairwise disjoint). We obtain reduced
$G$-trees $T_0, T'_0$ with the same elliptic subgroups: those
of $T$, and subgroups of an $N(G_v)$. Furthermore, $T_0$
satisfies the hypothesis of Theorem 1, with only cases
$(1)$ and $(2)$ occuring. Thus $T_0=T'_0$.

To reconstruct $T$ and $T'$ from $T_0$ and $T'_0$, one has to
blow up certain points into lines. Since these points project
onto terminal vertices in the quotient graph $T_0/G$, there is
only one way of blowing up. This shows $T'=T$.

\head Deformations \endhead

Given a $G$-tree $T$, a {\it collapse move\/} consists in
choosing an edge $e=vu$ such that $v,u$  are in distinct orbits
and $G_e=G_u$, and collapsing each edge in the  orbit of $e$ to
a point. The stabilizer of $v$ has not changed since
$G_v*_{G_u}G_u=G_v$. In the quotient graph of groups, one has
collapsed an edge to  a point. An {\it expansion move\/} is the
opposite operation.

Now we consider {\it slide moves\/} (see [\Fo] or [\FoJSJ]).
Suppose that  
$T$ contains  adjacent edges $e,f$ with origin $v$ such that
$G_e\inc G_f$, and $e,f$ are not in the same $G$-orbit as
non-oriented edges. We can then   slide $e$ across
$f$, so that it now starts at the terminal endpoint of $f$. Doing
this  
$G$-equivariantly replaces $T$ by another $G$-tree $T_1$.

\thm{Lemma} Let $T$ be reduced. If
$e,f$ are not as in case $(3)$ of Theorem 1, then sliding $e$ across
$f$  does change $T$ (the new tree $T_1$ is not equivariantly
isomorphic to
$T$).
\fthm

\demo{Proof} We may assume  that $T_1$ is reduced. We
show that the translation length of some element of $G$
changes.

Let
$u$ (resp\.
$w$) be the terminal endpoint of
$e$ (resp\.
$f$). If there is $g\in G$ sending
$v$ to $u$, the result is clear because the translation length of
$g$ goes from
$1$ to
$2$. Suppose there is no such $g$. In particular,  
$G_e\subsetneqq G_u$ since $T$ is reduced. Fix
$g_u\in G_u\setminus G_e$.

If $G_f\subsetneqq G_v$, choose   $h\in G_v\setminus G_f$.
Note that $h\notin G_e$ (because $G_e\inc G_f$). The
translation length of
$g_uh$ is twice the distance between the fixed point sets of
$g_u$ and $h$, so passes from
$2$   to
$4$. The case
$G_f\subsetneqq G_w$ is similar, so we assume
$G_f=G_v=G_w$.

If we are not in case $(3)$ of Theorem 1, there is an edge
$e_0=vv_0$ not in the same
$G$-orbit as
$e$ or $f$ (as a non-oriented edge). If there exists $g_0\in
G_{v_0}\setminus G_{e_0}$, the translation length of
$g_0g_u$ goes from $4$ to $6$. If
$G_{e_0}=G_{v_0}$, choose $h$ with $hv_0=v$. The translation
length of $g_uh$ goes from to
$3$ to $5$ (the translation length of $h$ is $1$, and the distance
from its axis to the fixed point set of $g_u$ goes from $1$ to
$2$).
\cqfd
\enddemo

We now prove the ``only if'' direction of Theorem 1, by
deforming $T$ into a
     reduced tree $T'$ different from $T$. Consider adjacent
edges $e,f$, with $G_e\inc G_f$,  satisfying none of the three
conditions of Theorem 1. There are two cases.

If
$e$ and $\bar f$ are not in the same orbit, we can change $T$
within its deformation space by sliding $e$ across $f$. The new
tree is reduced, except  if
$G_e=G_f=G_w\subsetneqq G_v$, where $w$ is the terminal
endpoint of $f$.   If this happens, we choose $t\in G$ taking 
$w$ to $v$ and we slide $e$ across   $t\bar f$ rather than across
$f$.

The second possibility is that $e$ and $\bar f$ are in the same
orbit and
$G_e\subsetneqq G_f$.  We may assume $G_f\subsetneqq G_v$:
otherwise
 $T$ is an ascending HNN-extension (if every edge with origin
$v$ is in the $G_v$-orbit of $e$ or $f$),  or there exist edges as
in the previous case (if there is a third orbit).   
 
\midinsert 
\centerline 
{\includegraphics[scale=.8]{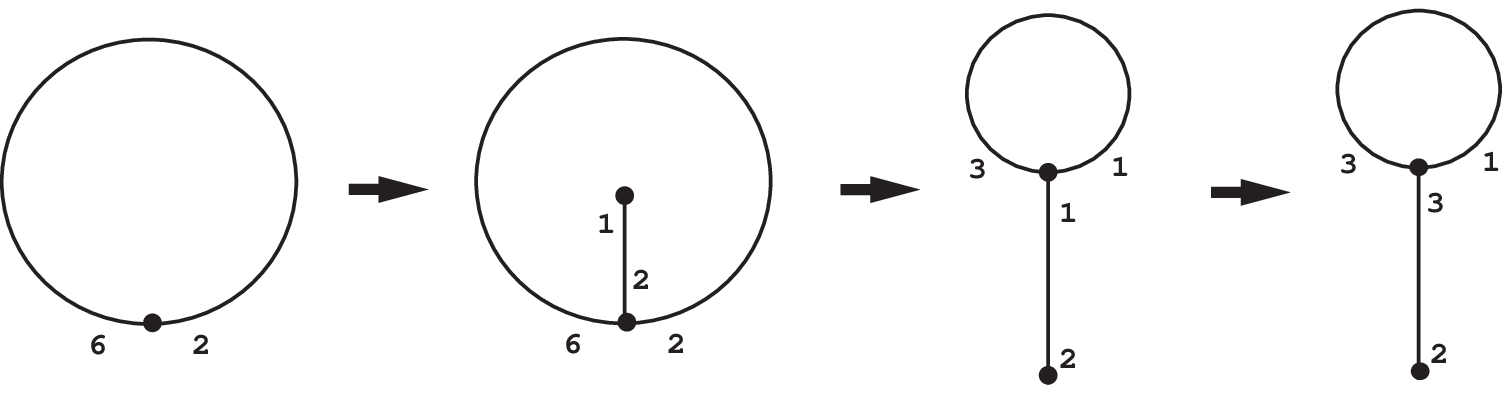}}
\botcaption
 {Figure 2}{Expanding and sliding in $BS(2,6)$}
\endcaption
\endinsert



We first perform an expansion move at $v$,   creating a new
edge with origin
$v$ and stabilizer $G_f$. In the quotient graph of groups, we
have created an edge and a terminal vertex, both with group
$G_f$. We then slide $e$ and $f$ across the new edge (the case
of $BS(2,6)$  is illustrated on figure 2). In the graph of groups,
the terminal vertex now has valence
$3$, and exactly two of the three edge groups map onto $G_f$
(the tree is not reduced). The last step is to   slide the new
edge across
$f$ (counterclockwise on figure 2). This yields a reduced tree
because $G_e\subsetneqq  G_f\subsetneqq G_v$.

\head Ascending HNN-extensions \endhead

We shall say that a $G$-tree $T_1$ is a {\it reduction\/} of a 
$G$-tree $T_2$ if it is reduced and may be obtained from
$T_2$ by performing collapse  moves. We also say that the
quotient graph of  groups $T_1/G$ is a reduction of $T_2/G$.
Every cocompact $G$-tree (resp\. every finite graph of
groups) has at least one reduction. 

  Consider a  graph of groups $\Gamma $ of the following form. 
It is a subdivided  circle, consisting of two vertices
$v,v'$ and two edges $e,e'$, and both inclusions
$G_{e}\hookrightarrow G_{v}$ and $G_{e'}\hookrightarrow
G_{v'}$ are onto. Such a graph of groups has two reductions,
obtained by collapsing  either edge. Both are ascending
HNN-extensions, and we say that the associated  Bass-Serre
trees are related by an {\it induction move\/} (through
$\Gamma $).

\thm{Lemma} Let $T$ be an ascending HNN-extension. If $T'$ is
a  reduced  tree in the deformation space of $T$, then $T'$ is an
ascending HNN-extension  and it may be obtained from
$T$ by a finite number of induction moves.
\fthm

\demo{Proof}  Join $T'$ to $T$ by a sequence $S_n$, with
$S_{n+1}$ obtained from $S_n$ by   an expansion or a collapse
move.  For each $n$, choose a reduction $T_n$ of  $S_n$. Thus
$T_n$ and $T_{n+1}$ are two different reductions of the same
tree ($S_n$ or
$S_{n+1}$). We complete the proof by showing that two
reductions of a tree $S$ in the deformation space of $T$ are
related by  an induction move (and are ascending
HNN-extensions).

Since $S$ may be obtained from $T$ by expansions and
collapses,   one can check that
   the graph of groups $\Gamma _S$ associated to   $S$ has the
following form. It consists of a (subdivided) circle $C$, possibly
with  trees $C_i$ attached to   vertices
$v_i\in C$. The circle $C$ may be oriented so that, if $e=vw$ is
an oriented  edge of
$C$, then $G_e=G_v$. Furthermore, the fundamental group of
$C$ (viewed as  a subgraph of groups) equals $G$ (i.e\. the
fundamental group of $C_i$ equals
$G_{v_i}$).

It follows from this description that any reduction of $S$ is an
ascending HNN-extension, obtained by choosing an edge
$e\inc C$ and collapsing  all other edges. Now consider two
reductions of
$S$. They are associated to edges $e,e'\inc C$, and they differ
by an  induction move (through the graph of groups $\Gamma $
obtained by collapsing all edges of
$\Gamma _S$ except
$e$ and $e'$).
\cqfd\enddemo

\thm{Corollary} $T$ is rigid if and only if $T'=T$ for every $T'$
related to
$T$ by an induction move. \cqfd
\fthm

  Let $T$ be associated to the
   ascending HNN-extension
$G =\langle A,t\mid tat\mi=\varphi (a)\rangle
$, with $\varphi :A\to A$ an injective homomorphism. It
contains an  edge $e=vw$, with $G_v=A$, $G_e=G_w=\varphi
(A)$, and $w=tv$.

Let $T'$ be related to $T$ by an induction move through
$\Gamma $. Let $T_0$ be the  Bass-Serre tree of
$\Gamma $. In $T_0$, the segment between $v$ and $w$ consists
of two  edges $vv_0$ and $v_0w$. The stabilizer of
$v_0w$ is $\varphi (A)$, and the  stabilizer of $vv_0$ is a group
$H$ with $\varphi (A)\inc H\inc A$. The tree $T'$  (obtained from
$T_0$ by collapsing edges in the orbit of
$v_0w$) is the Bass-Serre tree
$T_H$ associated to the presentation of $G$ as the ``induced'' 
HNN-extension $G =\langle H,t\mid tht\mi=\varphi_{|H}
(h)\rangle
$. Conversely, if $\varphi (A)\inc H\inc A$, the tree $T_H$ is 
related to $T$ by an induction move.

This shows that {\it $T'$ is related to $T$ by an induction move
if and only if $T'$ is a $T_H$, with $\varphi (A)\inc H\inc A$.\/} 
Proving Theorem 2 now reduces to showing that  $T_H=T$ if and
only if there exist $i,j,a_0,a_1$ as in  the statement of the
theorem.

The tree $T_H$ is characterized (up to $G$-equivariant
isomorphism)  by the existence of an edge $e'=v'w'$ with 
$G_{v'}=H$ and
$G_{e'}=G_{w'} $, and an element $t'\in G$ sending $v'$ to $w'$
such that
   $t'ht'{}\mi=\varphi (h)$ for $h\in H$.

If $T_H=T$, view $v'w'$ as an edge of $T$ and  fix
$g\in G$ taking the ``base edge''
$vw$ to
$v'w'$. Recall that $G_v=A$ and $w=tv$. The elements taking
$v'$ to $w'$ are those of the form
$ gta_1g\mi$, with
$a_1\in A$. Thus $T_H=T$ if and only if there exist $g\in G$ and 
$a_1\in A$ such that 
$$\left\{\aligned &gAg\mi=H   \\ &(gta_1g\mi ) 
gag\mi(gta_1g\mi)\mi=\varphi  (gag\mi) \quad \text {for
$a\in A$.}\endaligned\right .$$ 

Any $g\in G$ may be written $g=t^{-j}a_0t^i$, with $i,j\ge0$ and 
$a_0\in A$. Since
$t^iAt^{-i}=\varphi ^i(A)$ and $t^jHt^{-j}=\varphi ^j(H)$, one has 
$gAg\mi=H$ if and only if $a_0\varphi ^i(A)a_0\mi=\varphi
^j(H)$.

The other equation then becomes
  $$t^{-j}a_0t^i t a_1aa_1\mi t\mi t^{-i}a_0\mi t^j=\varphi
(t^{-j}a_0t^iat^{-i}a_0\mi t^j).$$ We rewrite it as  $$a_0\varphi
^{i+1}(a_1aa_1\mi)a_0\mi=t^j \varphi (t^{-j}a_0\varphi
^i(a)a_0\mi t^j)t^{-j}.
$$

Since $a_0\varphi ^i(a)a_0\mi \in\varphi ^j(A)$, the right-hand
side equals
$\varphi (a_0\varphi ^i(a)a_0\mi ) $, and we get   $$a_0\varphi 
^{i+1}(a_1)\varphi ^{i+1}(a)\varphi ^{i+1}(a_1\mi)a_0\mi=
\varphi (a_0)\varphi ^{i+1}(a)\varphi (a_0)\mi .$$ This
expresses that $\varphi (a_0\mi)a_0\varphi ^{i+1}(a_1)$
centralizes $\varphi ^{i+1}(A)$,  and concludes the proof of
Theorem 2.

\Refs 
\widestnumber\no{99}
\refno=0

\bref \by B. Bowditch \paper Cut points and canonical splittings
of hyperbolic groups\jour Acta Math.
\vol180\yr1998\pages145--186\endref

\bref  \by M. Forester\paper Splittings of generalized
Baumslag-Solitar  groups\jour Preprint
\endref

\bref   \by M. Forester\paper Deformation and rigidity of
simplicial group actions on trees\jour Geom.
\& Topol.\vol 6\yr2002\pages 219--267\endref

\bref   \by M. Forester\paper On uniqueness of JSJ
decompositions of finitely generated groups\jour Comm. Math.
Helv. \vol78\yr2003\pages 740--751 \endref

\bref  \by N.D. Gilbert, J. Howie, V. Metaftsis, E. Raptis\paper
Tree actions of automorphism groups\jour J. Group
Theory\yr2000\vol3\pages213--223
\endref

\bref   \by V. Guirardel\paper A very short proof of Forester's
rigidity result\jour Geom. \& Topol.\vol 7\yr2003\pages
321--328\endref

\bref\by V. Guirardel, G. Levitt\jour In preparation \endref

\bref\by G. Levitt \paper Automorphisms of   hyperbolic
groups and graphs of groups \jour Geom.
Dedic. (to appear)
\endref 

\bref\by  G. Levitt\paper On automorphisms of generalized
Baumslag-Solitar groups
\jour In preparation
\endref

\bref\by K. Whyte\paper The large scale geometry of the higher
Baumslag-Solitar groups\jour Geom. Funct. Anal. \yr
2001\vol11\pages 1327--1343\endref

\endRefs

\address   LMNO, umr cnrs 6139, BP 5186, Universit\'e de Caen,
14032 Caen Cedex, France.\endaddress\email 
levitt\@math.unicaen.fr{}{}{}{}{}\endemail

\enddocument